\documentclass[11pt]{amsart}

\usepackage{amsmath, amsthm, amssymb}
\usepackage{graphicx, array}
\usepackage{longtable}
\usepackage{epstopdf}
\usepackage{stmaryrd}
\usepackage{float}
\usepackage{hyperref}
\usepackage[a4paper, margin=1.2in]{geometry}
\usepackage[normalem]{ulem}

\usepackage{comment}
\usepackage{cancel}
\usepackage{soul}

\usepackage{color}

\definecolor{lggrr}{rgb}{0.5,0.5,0.5}

\newcommand{\aut}{\operatorname{Aut}}
\newcommand{\id}{\operatorname{id}}
\newcommand{\dom}{\operatorname{dom}}
\newcommand{\ran}{\operatorname{ran}}

\newcommand{\Lrel}{\mathrel{\mathcal{L}}}
\newcommand{\Rrel}{\mathrel{\mathcal{R}}}
\newcommand{\Hrel}{\mathrel{\mathcal{H}}}
\newcommand{\Drel}{\mathrel{\mathcal{D}}}

\newcommand{\paut}[1]{\operatorname{PAut}(#1)}
\newcommand{\is}[1]{\operatorname{PSym}(#1)}
\newcommand{\Sgrp}{\mathcal{S}}

\newtheorem{Thm}{Theorem}[section]
\newtheorem{Prop}[Thm]{Proposition}
\newtheorem{Lem}[Thm]{Lemma}
\newtheorem{Cor}[Thm]{Corollary}

\theoremstyle{definition}

\newtheorem{Rem}[Thm]{Remark}

\newtheorem{Probl}[Thm]{Problem}

\newenvironment{Proof}{\rm \trivlist\item[\hskip \labelsep{\bf
Proof.\quad}]}{\hfill\qed\par\medskip\endtrivlist}

\title{Inverse monoids of partial graph automorphisms}

\author{Robert Jajcay}
\email{robert.jajcay@fmph.uniba.sk}
\address{Comenius University, Bratislava, Slovakia, and University of Primorska, Koper, Slovenia
}
\thanks{The first author was partially supported by VEGA 1/0474/15, VEGA 1/0596/17,  VEGA 1/0423/20, APVV-15-0220, and by the Slovenian Research Agency (research projects N1-0038, N1-0062, J1-9108).}

\author{Tatiana Jajcayov\' a}
\email{tatiana.jajcayova@fmph.uniba.sk}
\address{Comenius University, Bratislava, Slovakia}
\thanks{The second author was partially supported by VEGA 1/0039/17 and VEGA 1/0719/18.}

\author{N\'ora Szak\'acs }
\email{szakacsn@math.u-szeged.hu}
\address{Bolyai Institute,
University of Szeged\\
Aradi v\'ertan\'uk tere 1.
H-6720 Szeged, Hungary}
\thanks{The third and fourth authors were partially supported by the National Research, Development and Innovation Office, grants no. K104251 and K115518.}

\author{M\'aria B. Szendrei }
\email{m.szendrei@math.u-szeged.hu}
\address{Bolyai Institute,
University of Szeged\\
Aradi v\'ertan\'uk tere 1.
H-6720 Szeged, Hungary}

\begin{document}

\begin{abstract}
A partial automorphism of a finite graph is an isomorphism between its vertex-induced subgraphs. 
The set of all partial automorphisms of a given finite graph forms an inverse monoid under composition {(of partial maps)}. 
We describe the algebraic structure of such inverse monoids by the means of standard tools of inverse semigroup theory,
{namely} Green's relations and some properties of the natural partial order, and give a characterization of inverse monoids which arise as inverse monoids of partial graph automorphisms. 
We extend our results to digraphs and edge-colored digraphs as well.
\end{abstract}

\maketitle

\section{Introduction}

The study of almost any type of combinatorial structures includes problems of distinguishing different objects via
the use of isomorphisms, and, more specifically, leads to questions concerning {the} automorphism groups of the
structures under consideration. 
The problem of determining the automorphism group for a specific combinatorial structure
or a class of combinatorial structures is a notoriously hard computational task whose exact complexity is
the subject of intense research efforts. 
For example, the time complexity of determining the automorphism group of a finite
graph (as a function of its order) has been the focus of a series of recent highly publicized lectures of
Babai \cite{Bab}.

Knowledge of the automorphism group of a specific combinatorial object allows one to make various claims
about its structure. 
As an extreme example, the fact that the
induced action of the automorphism group of a graph on the
set of 
unordered pairs of its vertices is transitive yields 
immediately that the graph under consideration is either the complete
graph $K_n$ or its complement $\widetilde{K_n}$. 
Similarly, vertex-transitivity of the
action of the automorphism group of a graph 
yields that each
vertex of the graph is contained in the same number of cycles
of prescribed length (e.g., \cite{Fil&Jaj}).

However, the usefulness of knowing the automorphism group
of a graph decreases with an increasing number of orbits of 
its action on the vertices. 
Revoking once again an extreme case, 
knowing that the automorphism group of a graph is trivial
yields almost no information about the structure of the graph 
(those who would like to disagree should consider the fact that 
almost all finite graphs have a trivial automorphism group, \cite[Corollary 2.3.3]{God&Roy}).

This observation appears to suggest rather limited use of 
algebraic tools in general graph theory. 
To counter this view, we propose
an approach that relies on the use of an algebraic theory often 
viewed as a generalization of group theory which applies to 
all finite graphs and relies on the concept of partial isomorphisms of combinatorial
structures \cite{Lau3,Lawson}.
Namely, we propose to study {\emph{the partial automorphism monoid $\paut\Gamma$ of a 
finite graph $\Gamma$}, which consists of all isomorphisms between the vertex-induced 
subgraphs of $\Gamma$, 
called partial automorphisms of $\Gamma$,} 
with the operation of {the usual} composition {of partial maps}.
{Unlike the case of automorphism groups, no finite graph of at least two vertices admits a trivial inverse monoid of partial 
automorphisms.
Indeed, partial identical maps are always partial automorphisms, and these already account for exponentially more elements 
than the number of vertices of the graph. In addition, as we will see, there are usually many more partial automorphisms.

In this paper, we make steps towards developing a general theory of 
partial automorphism monoids of combinatorial structures, in particular, 
of finite graphs, digraphs and edge-colored digraphs. 
Our main aim is to understand which inverse monoids arise as partial automorphism monoids of these structures.

We begin by reviewing {the} relevant concepts of graph theory in Section \ref{sec:graph},
{followed by the basics needed} from inverse semigroup theory in Section
\ref{sec:sem}. 
In Section \ref{sec:struc}, we {formulate the most important properties which determine} the structure of the partial automorphism monoid of a graph.
We then address two closely related classification problems concerning the 
partial automorphism monoids of graphs, digraphs
and edge-colored digraphs.

In Theorems \ref{thm:pautgraph} and \ref{thm:pautdigraph}, we describe those inverse monoids which arise as partial automorphism monoids of graphs, that is,
we characterize those inverse submonoids $\Sgrp$ 
of the so-called symmetric inverse monoid on $X$ 
for which a 
{graph (resp. digraph or edge-colored digraph) $\Gamma$ exists on the set of vertices $X$
such that $\paut\Gamma$ is {\em equal} to $\Sgrp$. 

This problem is similar to the more specialized problem
of the classification of finite groups that admit a Graphical
Regular Representation (for short, GRR). A finite group $G$ 
admits a GRR if there exists a
set of edges $E$ on $G$ with the property that the automorphism group 
of the graph $(G,E)$ 
is equal to $G_L$, the left multiplication regular representation
of $G$ in its action on itself \cite{God}. 
Cayley's left multiplication representation of groups has a direct analogue 
in the inverse semigroup theory in the form of the Wagner--Preston representation \cite[Theorem 1.5.1]{Lawson} which represents $\Sgrp$ as a set of partial permutations on itself. These partial permutations of $\Sgrp$ give rise to color-preserving partial automorphisms of its Cayley graph, just like in the group case, however, the Cayley graph of an inverse semigroup $\Sgrp$ has many more partial automorphisms than just those
induced by $\Sgrp$. 
Thus, an exact equivalent to the GRR 
classification would be the classification of the finite inverse monoids
$\Sgrp$ that admit a set of edges $E$ on $\Sgrp$,
{such that} the partial automorphism monoid of the graph $(\Sgrp,E)$
is equal to the partial permutation representation of $\Sgrp$ on itself 
given by the Wagner--Preston theorem. 
However,
by the observation made above concerning the number of vertices and the 
number of partial automorphisms of a graph, 
the existence of such a graph $(\Sgrp,E)$ is impossible for any non-trivial inverse monoid $\Sgrp$.


}

In Section \ref{sec:frucht2},
building on Theorems \ref{thm:pautgraph} and \ref{thm:pautdigraph}, we give a similar description in Theorems \ref{thm:pautgraph2} and \ref{thm:pautdigraph2} for the finite inverse monoids which are 
{\em isomorphic} to {partial automorphism monoids of} 
finite graphs, digraphs or edge-colored digraphs.
It turns out that the class of finite inverse monoids arising as {partial automorphism monoids of} (edge-colored di)graphs is very restrictive.
This is in stark contrast to the well-known result of Frucht \cite{fru} who proved that 
every finite group is isomorphic to the automorphism group of
a finite graph. 
A kind of `extended' Frucht theorem for finite inverse monoids has been
obtained in \cite{Nem}, where it is proved
that every finite inverse monoid is isomorphic to the partial weighted automorphism monoid 
of a finite weighted graph. The construction takes a cleverly modified version of the Cayley color graph of an inverse monoid $\Sgrp$, and adds weights to the vertices in order to avoid the partial automorphisms not arising from the Wagner--Preston representation. The resulting structure therefore has a partial automorphism monoid isomorphic to $\Sgrp$.

A slightly different approach is taken in \cite{Sie} in order to achieve a theorem of a similar flavor: here, the Cayley color graph of $\Sgrp$ remains unchanged, but the notion of 
partial automorphisms is restricted sufficiently so that only the partial automorphisms arising from the Wagner--Preston representation satisfy the definition. The inverse monoid of such partial automorphisms is therefore again isomorphic to $\Sgrp$.

We feel obliged to address
the similarities and differences between our paper and the recent
paper \cite{Chi&Ple} that appeared while we were preparing our
article. The authors of \cite{Chi&Ple} consider two types of inverse monoids
associated with a finite undirected graph which is allowed to have multiple edges and loops:
the first type being the 
inverse monoids of all partial automorphisms between any two
(not necessarily 
vertex-induced) 
subgraphs of such graphs, and the second
type being the inverse monoids considered in our paper.
The majority of the results obtained in \cite{Chi&Ple} deal with the 
partial automorphism monoids of the first type (not considered in our paper), 
and focus on the structure of their idempotents and ideals. 
None of the results obtained in our paper is included among the results 
contained in \cite{Chi&Ple}.

{Finally, let us mention papers devoted to finite monoids and semigroups and their relation to graphs that take a different approach from ours. Graphs play a fundamental role in studying certain classes of finite monoids and semigroups (see, e.g., \cite{Cam2}, \cite{Pin}, \cite{RhodesSteinberg}), while 
monoids attached to graphs have also been considered (see, e.g., \cite{Cam1}).
The monoids and semigroups appearing in this type of research consist almost
exclusively of transformations (e.g., transition monoids of automata, transformation semigroups generated by idempotent transformations coming from graphs, endomorphism monoids of graphs) rather than of partial permutations considered
in our paper.}

\section{Preliminaries from Graph Theory}\label{sec:graph}

A \emph{graph} is an ordered pair $\Gamma=(V,E)$, where $V$ is the set of vertices, and $E$ is the set of 
(undirected) edges, which is a set of $2$-element subsets of $V$. 
Similarly, a \emph{digraph} is an ordered pair $\Gamma=(V,E)$, with $V$ being the set of vertices, and $E$ the set of (directed) edges, consisting of ordered pairs of vertices. 
Graphs can naturally be regarded as special digraphs {where each edge $\{u,v\}$ of the graph is
replaced by} two directed edges of opposite directions, $(u, v)$ and $(v, u)$. 
Notice that digraphs can also have \emph{loops}, while graphs or digraphs associated to graphs cannot.
In both cases, we use the notation $V(\Gamma)$ and $E(\Gamma)$ for $V$ and $E$, respectively.
Furthermore, we also consider \emph{edge-colored digraphs}, that is, structures $\Gamma=(V, E_1, \ldots, E_l)$, where $\{1,\ldots, l\}$ is the set of colors, and $E_c\subseteq V \times V \ (c\in \{1,2,\ldots,l\})$ are the pairwise disjoint sets of (directed) 
edges of color $c$.
In this case, $V(\Gamma)$ and $E(\Gamma)$ stand for $V$ and $\bigcup_{c=1}^l E_c$, respectively.
For example, Cayley color graphs
constitute a natural class of edge-colored digraphs.

A \emph{partial automorphism} of an edge-colored digraph $\Gamma$ (as well as of a digraph or a graph) 
is an isomorphism between two vertex-induced subgraphs of $\Gamma$, that is, 
a bijection $\varphi\colon V_1 \to V_2$ between two sets of vertices $V_1, V_2 \subseteq V(\Gamma)$ such that 
any pair of vertices $u, v \in V_1$ satisfies the condition 
$(u,v) \in E_c$ if and only if $(\varphi(u), \varphi(v)) \in E_c$ for any color $c$. 
The set of {all} partial automorphisms of $\Gamma$ together with the operation of {the usual 
composition of partial maps} form an inverse monoid, which we denote by $\paut\Gamma$.
{This composition, as well as further concepts needed later and several basic properties of $\paut\Gamma$
are discussed} in detail in the next section.
In the rest of the paper, vertex-induced subgraphs are simply called induced subgraphs.

As argued in the introduction, we believe that developing a full-fledged theory of partial automorphism monoids of graphs might lead to discoveries in areas
of graph theory that have traditionally resisted the use of algebraic methods. As an example of one such possible impact, we mention the long-standing open problem 
called \emph{Graph Reconstruction Conjecture}, first introduced in \cite{Kel}. 
Given a finite graph $\Gamma=(\{v_1,\ldots,v_n\},E)$ of order $n$, let the \emph{deck of} $\Gamma$, 
$\operatorname{Deck}(\Gamma)$, be the multiset consisting of the subgraphs   
$\Gamma - v_i\ (1 \leq i \leq n)$ induced by ($n-1$)-vertex subsets of $\Gamma$. 
The Graph Reconstruction Conjecture predicts the unique `reconstructability' of any graph $\Gamma$ 
of order at least $3$ from its $\operatorname{Deck}(\Gamma)$, or, in other words, predicts that 
the multisets
$\operatorname{Deck}(\Gamma_1)$ and $\operatorname{Deck}(\Gamma_2)$, up to isomorphism, coincide
if and only if $\Gamma_1$ and $\Gamma_2$ are isomorphic
(in notation, $\Gamma_1 \cong \Gamma_2$)
for any pair of graphs $\Gamma_1,\Gamma_2$ of order at least $3$. 
This conjecture appears closely related to partial automorphisms \cite{Koc1,Ill}.
Namely, any two induced subgraphs $\Gamma - v_i$ and $\Gamma - v_j\ (i\neq j)$ admit at least one partial isomorphism 
$\varphi$ with domain of size $n-2$, namely the `identity' isomorphism between 
$(\Gamma - v_i)-v_j$ and $(\Gamma - v_j)-v_i$. 
Clearly, in the case where $\Gamma - v_i$ and $\Gamma - v_j$ admit exactly one partial isomorphism $\varphi$ 
with domain of size $n-2$, $\Gamma$ is reconstructable from $\Gamma - v_i$ and $\Gamma - v_j$ alone, 
by `gluing' together $v$ and $\varphi(v)$, for each $v$ in the domain of $\varphi$.  
Furthermore, two induced subgraphs $\Gamma - v_i$ and $\Gamma - v_j$ of 
$\Gamma$ (or sometimes the vertices $v_i, v_j$) are said to be \emph{pseudo-similar} if 
$\Gamma - v_i \cong \Gamma - v_j$,
but no automorphism of $\Gamma$ maps $v_i$ to $v_j$ and $\Gamma - v_i$ to $\Gamma - v_j$.
Or, using the language of partial automorphisms, $\Gamma - v_i$ and $\Gamma - v_j$ are pseudo-similar, 
if $\paut\Gamma$ contains a partial automorphism mapping 
$\Gamma - v_i$ to $\Gamma - v_j$ that cannot be extended to an automorphism of $\Gamma$. 
It has been claimed in \cite{lau} that the Graph Reconstruction Conjecture
holds for graphs containing no pseudo-similar vertices (more precisely, an alleged `proof' of the 
Graph Reconstruction Conjecture is claimed to have been submitted for publication 
that falsely assumed the non-existence of pseudo-similar vertices). 
This suggests that a valid proof of the Graph Reconstruction Conjecture might be found through better 
understanding of the finite graphs containing pseudo-similar vertices, or equivalently,
through understanding the graphs $\Gamma$ with $\paut\Gamma$ containing elements
with domain size $n-1$ that cannot be extended to automorphisms of $\Gamma$. 
While graphs whose decks split into pairs of pseudo-similar subgraphs can be easily
constructed from graphical regular representations of groups of odd order \cite{Kim&Sch&Sto},
graphs containing arbitrarily large subsets of mutually pseudo-similar subgraphs, while they exist 
\cite{Kim&Sch&Sto,Lau2}, 
are generally hard to find and tend to be large compared
to the size of their sets of mutually pseudo-similar subgraphs. 
The topic of constructing graphs containing a small number of mutually 
pseudo-similar vertices has also been the focus of intense research 
\cite{God&Koc1,God&Koc2,Lau4}. Inspired by these results, we suggest that the following problem
from \cite{Kim&Sch&Sto} may be accessible via the use of methods of  
partial automorphism monoids.

\begin{Probl} 
Determine the orders of the smallest graphs $\Gamma_k$ containing
$k \geq 2$ mutually pseudo-similar vertices.
\end{Probl}

\section{Preliminaries from semigroup theory}\label{sec:sem}

In this section, we provide an introduction to our main tool {needed in the paper}: the algebraic structure of {finite} inverse monoids. 
We follow \cite{Lawson}, a monograph addressed to a much wider readership than researchers interested in semigroup theory,
in its conventional notation of maps (functions):
given a map $\varphi$, the element assigned to an element $x$ is denoted $\varphi(x)$, and the composition of
maps $\varphi$ and $\psi$, where $\varphi$ is applied first and $\psi$ after that, is denoted $\psi\varphi$.
{The interested reader is referred to \cite{Lawson} for further details of the theory of inverse semigroups.}

\subsection{Partial permutations}

Given a set $X$, we call a bijection between two subsets $Y, Z \subseteq X$ a \emph{partial permutation} of $X$. 
{The set of all partial permutations of $X$ is denoted $\is X$.}
In the paper, we will always assume $X$ to be finite. We allow for $Y=Z=\emptyset$, in which case we obtain the empty map. If $\varphi\colon Y\to Z\in \is X$, then $Y$ and $Z$ are the \emph{domain} and \emph{range} of $\varphi$,
denoted by $\dom\varphi$ and $\ran\varphi$, respectively.
The common size $|\dom\varphi|=|\ran\varphi|$ of the sets $\dom\varphi$ and $\ran\varphi$ is called the \emph{rank} of $\varphi$.
{The inverse of $\varphi$, considered as a bijection, is also a partial permutation of $X$ denoted by $\varphi^{-1}$, and 
$\dom\varphi^{-1}=\ran\varphi$ and $\ran\varphi^{-1}=\dom\varphi$.}

Partial permutations can be composed: given two maps $\varphi_1\colon Y_1 \to Z_1$ and $\varphi_2\colon Y_2 \to Z_2$, one obtains their composition $\varphi_2 \varphi_1$ by composing them on the largest 
{subset of $X$}
where it `makes sense' to do so, that is, on 
${\dom \varphi_2\varphi_1=}\varphi_1^{-1}(Z_1 \cap Y_2)$, 
where by definition $(\varphi_2 \varphi_1)(x)=\varphi_2(\varphi_1(x))$ for any $x$. 
{The range of $\varphi_2 \varphi_1$ is
$\ran \varphi_2 \varphi_1=\varphi_2(Z_1 \cap Y_2)$.}
It may happen that $Z_1 \cap Y_2=\emptyset$, in which case $\varphi_2 \varphi_1$ is the empty map. 

One of the first theorems one learns about permutations is that they can be written as {compositions of pairwise disjoint cyclic permutations} in a unique way, where 
cyclic permutations correspond to the orbits.
An analogous theorem exists for partial permutations as well, as described in \cite{Lips}, although orbits are slightly more complicated,
and a partial permutation is decomposed into the union of its restrictions to the orbits rather than as a product.
We slightly deviate from the notation of \cite{Lips},
in particular, 
the right-left composition of maps followed in the paper induces notation which is
quite unusual both in semigroup and permutation group theories.

If $\varphi$ is a partial permutation of $X$ and $x \in \dom \varphi$, then either there exists some $n \in \mathbb{N}$ such that $\varphi^{n}(x)=x$, in which case
 the orbit of $x$ is permuted cyclically by $\varphi$ (which is always the case if $\varphi$ is a permutation), or it may be that $\varphi^{n}(x)$ eventually ends up outside $\dom\varphi$, in which case the orbit forms what is called a path.

A cyclic permutation 
cyclically permutes all the elements it does not fix. In the realm of partial permutations, this notion is not as useful as that of partial permutations which permute all elements in their domain cyclically. We call these partial permutations \emph{cycles}. Of course, one can turn a nonidentical cyclic permutation into a cycle by omitting fixed points from its domain, and vice versa.
The cycle $x_1\mapsfrom x_k\mapsfrom \cdots\mapsfrom x_2\mapsfrom x_1$, 
where $x_1, x_2,\ldots, x_k$ are pairwise distinct and $k \geq 1$, is denoted by $(x_k\, \cdots\, x_2\, x_1)$.
A {\emph{path}} is a partial permutation of the form 
$x_k\mapsfrom \cdots\mapsfrom x_2\mapsfrom x_1$, denoted by $[x_k\, \cdots\, x_2\, x_1)$, where $x_1, x_2, \ldots x_k$ are all distinct and $k \geq 2$. 
Notice that $\dom [x_k\, \cdots\, x_2\,x_1)=\{x_1, x_2,\ldots, x_{k-1}\}$ and 
$\ran [x_k\, \cdots\, x_2\,x_1)=\{x_2, x_3,\ldots, x_k\}$. 

To understand how partial permutations decompose into cycles and paths, consider for instance the following partial permutations of the set 
$X=\{1,\ldots, 6\}$: 
{
\begin{eqnarray*}
\varphi\colon \{1,2,3,4\} \to \{1,2,4,5\}&\!\!\!\!\!,&
1 \mapsto 2,\ 2 \mapsto 1,\ 3 \mapsto 4,\ 4 \mapsto 5;\\
\psi\colon \{1,2,3,4,5,6\} \to \{1,2,3,4,5,6\}&\!\!\!\!\!,&
1 \mapsto 2,\ 2 \mapsto 1,\ 3 \mapsto 4,\ 4 \mapsto 5,\ 5 \mapsto 3,\ 6 \mapsto 6.
\end{eqnarray*}}
{Notice that $\varphi$} 
decomposes into the cycle $(2\,1)$ and the path $[5\,4\,3)$, 
{where the sets $\{1,2\},\ \{3,4,5\}$ are the orbits of $\varphi$. That is, $\varphi$ is the union of these two partial permutations, which we denote by 
$({2\,1}) \vee [5\,4\,3)$.
Similarly, $\psi=(2\,1)\vee (5\,4\,3)\vee (6)$, and it is a permutation of $X$. 
Notice that the last member $(6)$ of this union cannot be deleted because the partial permutation $(2\,1)\vee (5\,4\,3)$ has $\{1,2,3,4,5\}$ as its common domain and range rather than $X$.}

Two partial permutations $\varphi,\psi\in \is X$ are said to be {\em disjoint} if
the unions $\dom\varphi \cup \ran\varphi$ and $\dom\psi  \cup  \ran\psi$ are disjoint.
It is true that every partial permutation in $\is X$ is {the union of pairwise disjoint paths and cycles}, 
and the set of the members of this union is unique.
Using this form, one can  multiply partial permutations similarly as it is usually done in the case of permutations (from right to left)
{by starting with an element of the domain of the right factor, finding its image, then finding the image of this element, etc.,
but paths make the calculation somewhat more complicated because the choice of a starting element determines whether the whole path is obtained or only a part of it.}   
For example, 
$([4\,3\,1) \vee (2))([4\,1) \vee (3\,2))=[4\,2)\vee [4\,2\,3)=[4\,2\,3)$.

The inverse of a partial permutation is obtained by reversing all paths and cycles, e.g.\ 
$(({2\,1}) \vee [5\,4\,3))^{-1}=({1\,2})
\vee [3\,4\,5)$.

To avoid having to separate too many cases later in our proofs, we will sacrifice uniqueness and allow for the possibility of $[x_k\, \cdots\, x_2\,x_1)$ denoting a cycle with $x_1=x_k$, but still indicate cycles in the usual form 
when it is more convenient.
For example, this will allow us to write any $\varphi\in \is X$ with $\dom \varphi=\{x_1, x_2,\ldots,x_k\}$
in the form $\varphi=[\varphi(x_1)\, x_1)\vee [\varphi(x_2)\, x_2)\vee \cdots \vee [\varphi(x_k)\, x_k)$.

Let us further explore the partial operation $\vee$ we used in the decomposition. Given arbitrary partial permutations $\varphi$ and $\psi$, their union is a map if an only if they coincide on their common domain $\dom \varphi \cap \dom \psi$. In this case, this map is injective, that is, a partial permutation itself if and only if 
$\varphi$ and $\psi$ map elements outside their common domain into disjoint sets. In this case, we call $\varphi$ and $\psi$ \emph{compatible}, and denote their union by $\varphi \vee \psi$. This is in fact the least upper bound of $\varphi$ and $\psi$ in the partial order of restriction, defined by $\varphi_1 \leq \varphi_2$ if $\varphi_1$ is the restriction of $\varphi_2$ to some set $Y \subseteq X$, that is, if $\varphi_1=\varphi_2\id_Y$. For this reason, from now on we refer to $\varphi \vee \psi$ as the \emph{join} of $\varphi$ and $\psi$. Given an arbitrary subset $S$ of $\is X$, its join exists if and only if the elements of $S$ are pairwise compatible, and in this case, its join is the union of its elements.

The operations of composition and taking inverse {in $\is X$} distribute over joins. 
More precisely,
for any $\varphi,\psi,\eta\in\is X$,
if $\varphi\vee\psi$ exists then all of 
$\varphi\eta\vee\psi\eta$, $\eta\varphi\vee\eta\psi$, and $\varphi^{-1}\vee\psi^{-1}$ exist, 
and we have
$(\varphi\vee\psi)\eta=\varphi\eta\vee\psi\eta$,
$\eta(\varphi\vee\psi)=\eta\varphi\vee\eta\psi$, and
$(\varphi\vee\psi)^{-1}=\varphi^{-1}\vee\psi^{-1}$.

\subsection{Inverse monoids}\label{ssec:invmon}
 
The concept of inverse monoids was independently defined by Wagner and Preston in the early 50s as the algebraic abstraction of partial permutations. 
The set {$\is X$} of {all} partial permutations of a set $X$ forms a monoid under composition, $\id_X$ being the identity element. 
Elements also have `local' inverses: for any partial permutation 
$\varphi\colon Y \to Z$, its inverse $\varphi^{-1}\colon Z \to Y$ has the property 
that $\varphi^{-1}\varphi=\id_{Y}$ and $\varphi\varphi^{-1}=\id_{Z}$, and therefore 
the identities $\varphi\varphi^{-1}\varphi=\varphi$ and $\varphi^{-1}\varphi\varphi^{-1}=\varphi^{-1}$ hold. 
Moreover, $\varphi^{-1}$ is the only element with these properties.

Generalizing the above, a monoid $\Sgrp$ {is said to be} an 
\emph{inverse monoid} if for every $s \in \Sgrp$ there exists a unique element 
$s^{-1} \in \Sgrp$ called the inverse of $s$ such that 
$ss^{-1}s=s$ and $s^{-1}ss^{-1}=s^{-1}$ hold. 
Note that the unary operation of taking inverse has the properties $(s^{-1})^{-1}=s$ and $(st)^{-1}=t^{-1}s^{-1}$ for any $s,t\in \Sgrp$.
{The inverse monoid of all partial permutations of $X$ is called the \emph{symmetric inverse monoid on $X$}, and is also denoted by $\is X$.}

Symmetric inverse monoids are the archetypal examples of inverse monoids.
They are archetypal in the same sense as symmetric groups are the archetypal groups: in parallel with Cayley's theorem, the Wagner--Preston theorem states that every inverse monoid can be embedded into a suitable symmetric inverse monoid.

Another example of an inverse monoid is the set $\paut \Gamma$ of partial automorphisms of a finite edge-colored digraph (graph, digraph) $\Gamma$;
we again allow for the empty map as well. Each partial automorphism of $\Gamma$, that is, an isomorphism between two induced subgraphs of $\Gamma$, is a partial permutation of $V(\Gamma)$, thus, $\paut \Gamma \subseteq \is {V(\Gamma)}$. It is easy to see that $\paut \Gamma$ is closed under partial {composition} and taking inverses and contains the identity element $\id_{V(\Gamma)}$, thus, $\paut \Gamma$ 
forms an \emph{inverse submonoid} of $\is {V(\Gamma)}$.
If $\varphi\in \paut \Gamma$ then the subgraphs of $\Gamma$ induced on the vertex sets $\dom \varphi$ and $\ran \varphi$ are also denoted by $\dom \varphi$ and $\ran \varphi$, respectively.

The partial order discussed above  induced by restriction on $\is X$ can be extended to any inverse monoid $\Sgrp$ as follows. First, consider the set $E(\Sgrp)$ of \emph{idempotent} elements in $\Sgrp$ satisfying $e^2=e$. These are the elements arising in the form $ss^{-1}$ {for some $s\in \Sgrp$}. Idempotents in $E(\Sgrp)$ necessarily 
commute and the product of two idempotents is again an idempotent. 
{Thus, the set of idempotents $E(\Sgrp)$ forms a semilattice under multiplication.}
Furthermore, an idempotent is always its own inverse, which means that 
$E(\Sgrp)$ forms an inverse submonoid. 

In the case of $\Sgrp=\is X$, the set of idempotents is $E(\Sgrp)=\{\id_{Y}\colon Y \subseteq X\}$  with
$\id_\emptyset$ being the empty map.
If $Y, Z \subseteq X$, then $\id_Y \id_Z=\id_{Y \cap Z}$, so the subsemilattice {$E(\Sgrp)$ is} isomorphic to the semilattice $(\mathcal P(X), \cap)$. 
The semilattice structure defines a partial order on $E(\Sgrp)$, and this partial order can be naturally extended to $\Sgrp$ by putting $s \leq t$ if and only if there exists $e \in E(\Sgrp)$ such that $s=te$. This is called the \emph{natural partial order} on $\Sgrp$. Notice that this is indeed the restriction order on $\is X$. {Moreover, if} $\Gamma$ is a finite edge-colored digraph, then  $E(\is{V(\Gamma)})\subseteq \paut{\Gamma}$, therefore $E(\paut{\Gamma})=E(\is{V(\Gamma)})$. So the induced natural partial order on $\paut{\Gamma}$ is also the restriction order.

The concept of compatibility we introduced in $ \is X$ can also be extended to abstract inverse monoids. 
Two partial permutations $\varphi$ and $\psi$ in $\is X$ coincide on their common domain if and only if $\varphi^{-1}\psi$ is an idempotent, while $\varphi$ and $\psi$ map elements outside their common domain into disjoint sets if and only if $\varphi\psi^{-1}$ is an idempotent. In an abstract inverse monoid $\Sgrp$, we will call a pair of elements $a,b \in \Sgrp$ \emph{compatible} if $ab^{-1}$ and $a^{-1}b$ are idempotents. A \emph{subset} of $\Sgrp$ is called \emph{compatible} if all its elements are pairwise compatible. Although every compatible subset of a symmetric inverse monoid
has a join (a least upper bound), compatibility is not a sufficient condition for the
existence of a join in general inverse monoids. For instance, let $\Gamma_0$ be the graph in Figure \ref{fig:graph}, and consider the elements $[1\,2)$ and $[4\,3)$ of $\paut {\Gamma_0}$. These are compatible, as $[1\,2)[4\,3)^{-1}=[1\,2)^{-1}[4\,3)=\id_\emptyset$, but any common upper bound of $[1\,2)$ and $[4\,3)$ in the restriction order would need to take the edge $\{2,3\}$ to the non-edge $\{1,4\}$, 
and so their join in $\paut{\Gamma_0}$ does not exist.

\begin{figure}
\includegraphics[width=0.25\linewidth]{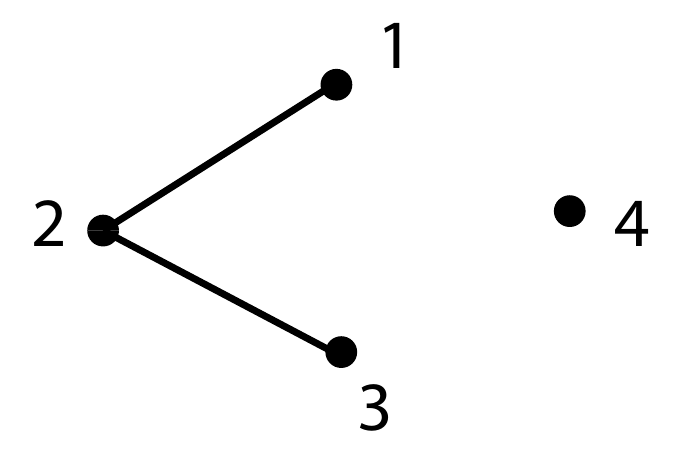}
\caption{The graph $\Gamma_0$}
\label{fig:graph}
\end{figure}

An inverse monoid $\Sgrp$ with zero is called {\em Boolean} if the semilattice $E(\Sgrp)$ is the meet semilattice of a Boolean algebra. 
{Notice that the zero element of a Boolean inverse monoid $\Sgrp$ is necessarily the minimum element of $E(\Sgrp)$.}
Both symmetric inverse monoids and partial automorphism monoids {of edge-colored digraphs (for brevity, partial automorphism monoids)} are Boolean: in both cases, the semilattice of idempotents is isomorphic to $(\mathcal P(X), \cap)$, and $\mathcal P(X)$ is of course a Boolean algebra. Furthermore, the minimum element $\id_{\emptyset}$ is indeed a zero element for the whole monoid. 

In this paper, we will mainly be concerned with finite inverse monoids: if $\Sgrp$ is finite, the fact that $E(\Sgrp)$ is Boolean automatically implies that $E(\Sgrp)$ is isomorphic to 
$(\mathcal P(X), \cap)$ for some finite set $X$.

We say that $\mathcal{T}$ is a \emph{full inverse submonoid} of $\Sgrp$ if $\mathcal{T}$ is 
an inverse submonoid of $\Sgrp$ such that 
$E(\mathcal{T})=E(\Sgrp)$. 
For instance, we have seen that $\paut{\Gamma}$ is a full inverse submonoid of $\is {V(\Gamma)}$. Some of our previous assertions were a direct consequence of this fact, for instance that the natural partial orders of $\paut{\Gamma}$ and $\is {V(\Gamma)}$ coincide on $\paut{\Gamma}$, or that $\paut{\Gamma}$ is Boolean. Notice that any full inverse submonoid of a symmetric inverse monoid $\is X$ is Boolean, and is closed under taking restrictions.

\subsection{Green's relations}\label{ssec:green}

A key to understanding the structure of inverse monoids (or semigroups in general) is knowing which elements can be multiplied to which elements. Some information is carried by the natural partial order -- one can never move `up' by multiplication in finite inverse monoids -- but an accurate description of this is best captured by five equivalence relations,
called \emph{Green's relations}, two of which coincide in the cases we are concerned with.

Let $\Sgrp$ be {an arbitrary monoid}, 
and $a, b \in \Sgrp$ be arbitrary elements. We define two equivalence relations $\Lrel$ and $\Rrel$ the following way:

\smallskip
\begin{center}
	$a \Lrel b$ if and only if there exist $x,y \in \Sgrp$ such that $xa=b$ and $yb=a$,\\
	$a \Rrel b$ if and only if there exist $x,y \in \Sgrp$ such that $ax=b$ and $by=a$.
\end{center}

\smallskip
\noindent
In a symmetric inverse monoid for instance, note that $\dom \varphi\psi \subseteq \dom \psi$ for any pair of partial permutations $\varphi$, $\psi$, so if $\varphi_1 \Lrel \varphi_2$, then clearly $\dom \varphi_1=\dom \varphi_2$ holds. On the other hand if $\dom \varphi_1=\dom \varphi_2$, then 
$(\varphi_2\varphi_1^{-1})\varphi_1=\varphi_2 \id_{\dom\varphi_1}=\varphi_2$, and similarly, $(\varphi_1\varphi_2^{-1})\varphi_2=\varphi_1 \id_{\dom\varphi_2}=\varphi_1$, so $\varphi_1 \Lrel \varphi_2$. Therefore, the $\Lrel$ relation coincides with having the same domain. Dually, the $\Rrel$ relation coincides with 
pairs of partial permutations having the same range. Note that the above argument holds without any modifications in $\paut \Gamma$ as well (or in any other inverse submonoid of 
a symmetric inverse monoid).
So, for instance, returning to $\paut {\Gamma_0}$ discussed above, the $\Rrel$-class of the partial automorphism $(1)\vee(2)$ consists of all partial automorphisms of range $\{1,2\}$, that is, of $\{(1)\vee(2), (1\,2), [1\,2\,3), (2) \vee [1\,3)\}$. Similarly, the $\Lrel$-class of $(1)\vee(2)$ consist of all partial automorphisms with domain $\{1,2\}$, that is, of $\{(1)\vee(2), (1\,2), [3\,2\,1), (2) \vee [3\,1)\}$. Note that the $\Lrel$-class of $(1)\vee(2)$ contains exactly the inverses of the elements in its $\Rrel$-class. This is no coincidence: in any inverse monoid, $a \Rrel b$ if and only if $a^{-1} \Lrel b^{-1}$.

The third Green's relation is ${\Hrel} = {\Rrel \cap \Lrel}$, which is again an equivalence relation. In a symmetric inverse monoid, or in a partial automorphism monoid, we clearly have $\varphi_1 \Hrel \varphi_2$ if and only if $\dom \varphi_1=\dom \varphi_2$ and $\ran \varphi_1 =\ran \varphi_2$. 
The $\Hrel$-class of $(1) \vee (2)$ in $\paut {\Gamma_0}$ is $\{(1)\vee(2), (1\,2)\}$. Note that it forms a subgroup. This is again not a coincidence.
In an inverse monoid, each $\Rrel$-class and each $\Lrel$-class are known to contain precisely one idempotent, and 
the $\Hrel$-classes containing these idempotents are the maximal subgroups of the inverse monoid. 
This is easy to see in 
partial automorphism monoids: if $\Gamma$ is an edge-colored digraph and $\Delta$ is an induced subgraph, then there is exactly one idempotent with a prescribed domain (or range) $V(\Delta)$, namely $\id_{V(\Delta)}$. The $\Hrel$-class of 
the idempotent $\id_{V(\Delta)}$ consists of all partial automorphisms $\varphi$ with $\dom\varphi=\ran\varphi=V(\Delta)$, which is exactly the automorphism group of the graph $\Delta$.

The smallest equivalence relation containing both $\Rrel$ and $\Lrel$ is called the $\Drel$ relation, and in any monoid, ${\Drel}={\Rrel \circ \Lrel} = {\Lrel \circ \Rrel}$, 
where $\circ$ denotes the composition of relations. 
Thus, for any $a,b \in \Sgrp$, we have $a \Drel b$ if and only if there exists $c \in \Sgrp$ such that {$a \Lrel c \Rrel b$}. In a symmetric inverse monoid, this means 
that $\varphi_1 \Drel \varphi_2$ if and only if there exists a partial permutation $\psi$ with $\dom\psi=\dom\varphi_1$ and $\ran\psi=\ran\varphi_2$, which is equivalent to
$|\dom\varphi_1|=|\ran\varphi_2|$, or, in other words, $\varphi_1$ and $\varphi_2$ having the same rank.

The $\Drel$ relation has a different characterization in partial automorphism monoids. It is still necessary to require that $|\dom\varphi_1|=|\ran\varphi_2|$, but it is no longer sufficient. For instance, consider $\Gamma_0$ and the partial automorphisms 
$(12), (34) \in \paut{\Gamma_0}$
again. Both have rank $2$, but there is no partial graph automorphism $\psi \in \paut{\Gamma_0}$ with $\dom\psi=\{1,2\}$ and $\ran\psi=\{3,4\}$, as $\{1,2\}$ {is an edge} 
while $\{3,4\}$ is not. It turns out that {in}
partial automorphism monoids, the $\Drel$ relation corresponds to isomorphism classes of {induced} subgraphs, as formulated in the following proposition.

\begin{Prop}
\label{prop:d}
For any edge-colored digraph $\Gamma$, the $\Drel$-classes of $\paut{\Gamma}$ correspond to the isomorphism classes of 
{the} induced subgraphs of\/ $\Gamma$, that is, two elements are $\Drel$-related if and only if the subgraphs induced 
{on their domains (or ranges)} are isomorphic.
\end{Prop}

\begin{Proof}
Let $\varphi_1, \varphi_2 \in \paut{\Gamma}$. Then $\varphi_1 \Drel \varphi_2$ if and only if there exists $\psi \in \paut{\Gamma}$ with $\varphi_1 \Lrel \psi \Rrel \varphi_2$.
The latter relation is equivalent to the equations $\dom{\varphi_1}=\dom{\psi}$ and $\ran{\varphi_2}=\ran{\psi}$. 
Thus, such $\psi$ exists in $\paut{\Gamma}$ if and only if there is a graph isomorphism between the induced subgraphs on $\dom{\varphi_1}$ and $\ran{\varphi_2}$, which proves our statement.
\end{Proof}

In an inverse monoid, each $\Drel$-class is a disjoint union of $\Rrel$-classes, and a disjoint union of $\Lrel$-classes, and any $\Rrel$- and $\Lrel$-class of a $\Drel$-class intersect in an $\Hrel$-class.
Moreover, every $\Drel$-class contains the same number of $\Rrel$-classes as $\Lrel$-classes.
A $\Drel$-class is therefore usually depicted in what is called an `eggbox' diagram:
the rows are the $\Rrel$-classes, the columns are the $\Lrel$-classes, the small rectangles are the $\Hrel$-classes, 
and they are arranged in such a way that the $\Hrel$-classes containing idempotents are on the main diagonal.

For an example, see Figure~\ref{fig:D} which depicts the $\Drel$-class of $\paut{\Gamma_0}$ corresponding to the edges of the graph. The first row (column) corresponds to the range (domain) $\{1,2\}$, and the second row (column) corresponds to the range (domain) $\{2,3\}$. These belong to the same $\Drel$-class since
both $\{1,2\}$ and $\{2,3\}$ represent edges of  $\paut{\Gamma_0}$, and therefore
there exists a partial automorphism of rank $2$ mapping $\{1,2\}$ to $\{2,3\}$. The idempotents are colored grey.

\begin{center}
\begin{figure}[H]
\includegraphics[width=0.3\linewidth]{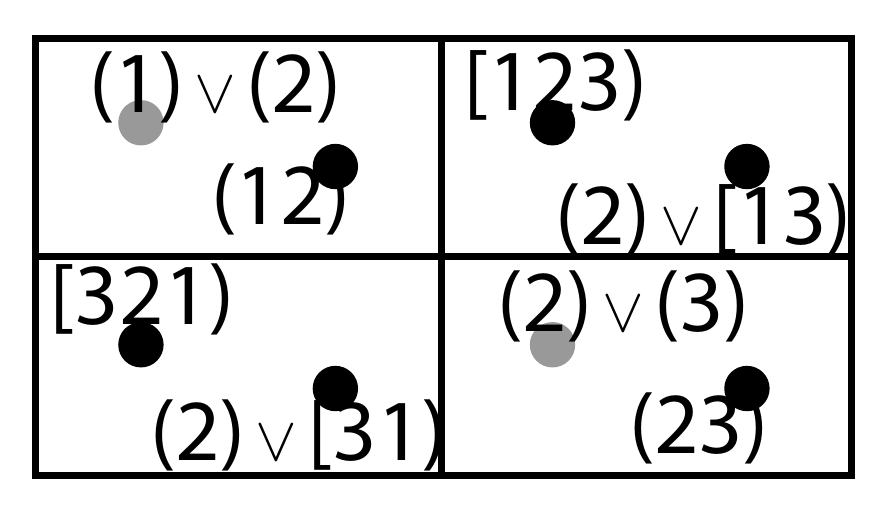}
\caption{The $\Drel$-class of {the} edges {in} $\paut{\Gamma_0}$}
\label{fig:D}
\end{figure}
\end{center}

In the case of finite semigroups, the $\Drel$-classes form a partially ordered set: if we 
denote the $\Drel$-class of $a$ by $D_a$, we put $D_a \leq D_b$ if $a=xby$ for some $x,y\in\Sgrp$, that is, if $b$ can be multiplied into $a$. 
In the case of $\is{X}$, this is just the ordering of $\Drel$-classes according to their rank; the $\Drel$-classes form a chain. 
The minimum element is the $\Drel$-class of partial permutations of rank $0$, which is just $\{\id_{\emptyset}\}$. 
The maximum element is the $\Drel$-class of rank {$|X|$} partial permutations, which are the permutations of $X$. 
This maximal $\Drel$-class therefore consists of a single $\Hrel$-class, which is the symmetric group $\operatorname{Sym}(X)$.

In the case of partial automorphism monoids, this partial order corresponds to the `induced subgraph' relation between the isomorphism classes of graphs,
as implied by the following proposition.

\begin{Prop}
Let $\Gamma$  be an edge-colored digraph, and let  $\varphi_1, \varphi_2 $
be elements of $ \paut{\Gamma}$. Then, $D_{\varphi_1} \le D_{\varphi_2}$ (i.e.,
there exist partial permutations $\psi, \sigma \in \paut{\Gamma}$ such that $\varphi_1 =\psi\varphi_2\sigma$)
if and only if $\dom\varphi_1$ is isomorphic to an induced subgraph of $\dom\varphi_2$.
\end{Prop}

\begin{Proof}
First suppose such $\psi$ and $\sigma$ exist. Then, 
$$\dom\varphi_1=\sigma^{-1}(\dom(\psi\varphi_2)\cap \ran \sigma)\subseteq\sigma^{-1}(\dom(\psi\varphi_2)) \subseteq \sigma^{-1}(\dom\varphi_2),$$ 
and the assertion holds.

Conversely, suppose that $\dom\varphi_1$ is isomorphic to an induced subgraph of $\dom\varphi_2$, {and} denote the {set of}
vertices of this subgraph by $W$. Then there exists a partial automorphism $\sigma \colon \dom \varphi_1 \to W$. Note that $\varphi_2(W)$ is isomorphic to $\ran\varphi_1$, so there exists a partial automorphism $\psi\colon \varphi_2(W) \to \ran\varphi_1$. Thus, $\varphi_1=\psi\varphi_2\sigma$.
\end{Proof}

In $\paut\Gamma$, the minimum element of the poset of $\Drel$-classes is again the singleton $\Drel$-class $\{\id_{\emptyset}\}$, while the maximum element is 
the single $\Hrel$-class equal to the automorphisms group $\aut(\Gamma)$.

If $\Sgrp$ is any finite semigroup with a zero element $0$, the $\Drel$-class $D_0=\{0\}$ 
is always the minimum element of this poset.
{I}n a {finite} poset $P$ with minimum element $0$, the {\em height} of an element $a$
is the largest 
	$s\in\mathbb{N}_0$ such that there exist elements $a_1,a_2,\ldots,a_s\in P$ with
	$0<a_1<a_2<\ldots<a_s=a$. An element of height $1$ is usually called \emph{$0$-minimal}.
This allows us to define the height of $\Drel$-classes in  $\Sgrp$, and also the height of elements in $\Sgrp$ with respect to the natural partial order. 
In any {finite}
inverse monoid with zero, the height of each element is known to be equal to the height of its $\Drel$-class. In $\is{X}$, the height of a partial permutation is exactly its rank, and the same is true in $\paut\Gamma$.

The $\Drel$-class depicted in Figure \ref{fig:D} is a $\Drel$-class of height $2$ in $\paut{\Gamma_{0}}$. There is one other $\Drel$-class of height $2$: the class that belongs to the isomorphism class {of} \emph{non-edges} (that is, pairs of vertices with no edge between them). Figure \ref{fig:green} illustrates the entire eggbox structure of $\paut{\Gamma_0}$.
The graph isomorphism class corresponding to a $\Drel$-class is depicted to the right of
it.

\begin{center}
\begin{figure}
\includegraphics[width=\linewidth]{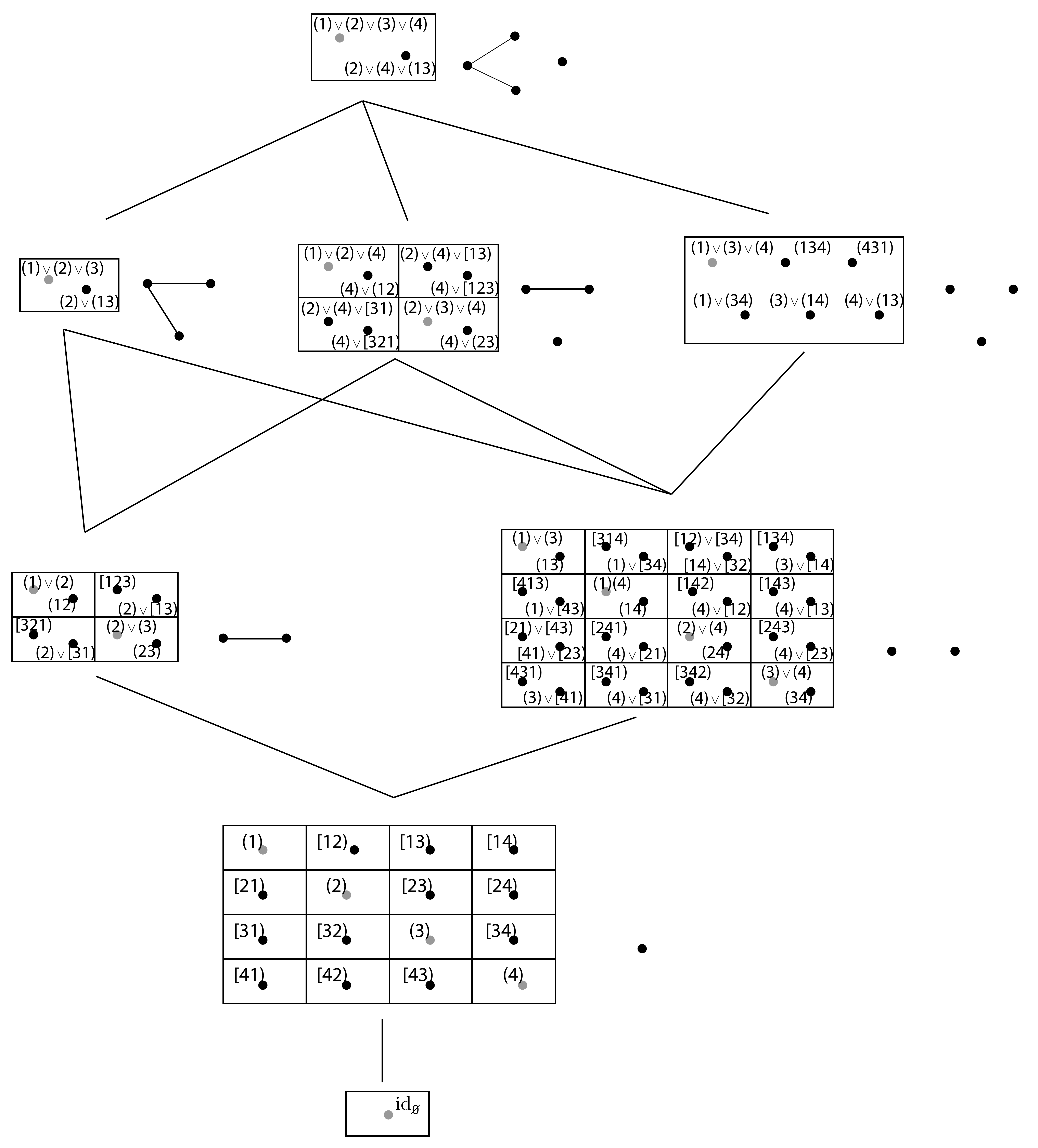}
\caption{The structure of the partial automorphism monoid $\paut{\Gamma_0}$}
\label{fig:green}
\end{figure}
\end{center}

\section{The structure of partial automorphism monoids of graphs} \label{sec:struc}

We begin by some simple observations. 
The partial automorphism monoids of a graph $\Gamma$ and its complement $\widetilde{\Gamma}$ are equal, that is,
$ \paut{\Gamma} = \paut{\widetilde{\Gamma}} $.
In particular, the partial automorphism monoid{s} of
the complete graph $K_n$ and its complement $\widetilde{K_n}$ are equal to
the symmetric inverse monoid $\is{V(K_n)}$. 
If the graphs $ \Gamma $ and $ \Gamma' $ are isomorphic, then $\paut{\Gamma}$ and $\paut{\Gamma'}$ are also
isomorphic (more specifically, $\paut{\Gamma} = \varphi^{-1} \paut{\Gamma'} \varphi $ for {any} isomorphism 
$ \varphi\colon \Gamma \to \Gamma' $).

A key observation in describing the partial automorphism monoids of edge-colored digraphs (in particular, of graphs) is the fact that the  elements of rank $1$ and $2$ determine the rest of the monoid. This is formalized in the following lemmas. The proof of the first lemma is obvious.

\begin{Lem}
\label{lem:rank2}
Le{t $\Gamma=(X,E)$ b}e an edge-colored {di}grap{h}, and  let $\varphi \in \is{X}$ be a partial permutation of rank at least $2$. 
Then $\varphi \in \paut{\Gamma}$ if and only if $\varphi|_{Y} \in \paut{\Gamma}$ for any $2$-element subse{t $Y$ of $\dom{\varphi}$.}
\end{Lem}

Lemma \ref{lem:rank2} implies that one can build all partial automorphisms of an edge-colored digraph (in particular, of a graph) from partial automorphisms of rank at most $2$, using joins.
\begin{Prop}
\label{prop:joins1}
The partial automorphism monoid $\Sgrp=\paut{\Gamma}$ of any edge-colored digraph $\Gamma$ has the following property: 
\begin{enumerate}
\item[{\rm (U)}]
for any compatible subset $A \subseteq \Sgrp$ of partial permutations of rank $1$, if $\Sgrp$ contains the join of any two elements of $A$, then $\Sgrp$ contains the join of the set $A$.
\end{enumerate}
\end{Prop}

\begin{Proof}
Since $A$ is compatible{, $\varphi=\bigvee A$ i}s an element of $\is{V(\Gamma)}$. 
The elements of $A$ have rank $1$, therefore compatibility of $A$ implies that distinct elements of $A$ have distinct domains and distinct ranges. 
Thus, for any distinct element{s $\psi_1, \psi_2 \in A$, $\psi_1 \vee \psi_2$ i}s a rank $2$ partial permutation 
{on} $V(\Gamma)$
(and, by assumption, also a partial graph automorphism of $\Gamma$), 
namely the restriction of $\varphi$ to the $2$-element set 
$\dom{\psi_1} \cup \dom{\psi_2}$. 
Moreover, by definition, any rank $2$ restriction of $\varphi$ arises in this way. 
Therefore $\varphi$ satisfies the conditions of Lemma \ref{lem:rank2}, and hence
belongs to $\Sgrp$.
\end{Proof}

\begin{Prop}
\label{prop:joins2}
If $\Sgrp,\mathcal{T}$ are full inverse submonoids of $\is{X}$ which coincide on their elements of rank at most $2$ and satisfy condition {\rm (U)}, then 
$\Sgrp =\mathcal{T} $. 
\end{Prop}

\begin{Proof}
To verify $\Sgrp \subseteq \mathcal{T}$, let $\varphi \in \Sgrp$ be of rank greater than $2$. 
Then all restrictions of $\varphi$ are in $\Sgrp$, since $\Sgrp$ is {a} full
{inverse submonoid of $\is{X}$}. 
In particular, with $A$ denoting the set of rank $1$ restrictions of $\varphi$, we have 
$A \subseteq \Sgrp$. 
The elements of $A$ are, of course, pairwise compatible, and as all their two-element joins are rank $2$ restrictions of $\varphi$, they belong to $\Sgrp$ as well. 
Since we assume that $\mathcal{T}$ contains the same elements of ranks $1$ and $2$
as $\Sgrp$ and that $\mathcal{T}$ satisfies condition {\rm (U)}, 
this yields $\varphi=\bigvee A $ is an element of $ \mathcal{T}$, as needed. 
The inclusion $\mathcal{T} \subseteq \Sgrp$ is obtained by the symmetric argument, swapping the roles of $\mathcal{T} $ and $ \Sgrp$.
\end{Proof}

\section{When does an inverse monoid of partial permutations coincide with the partial automorphism monoid of a graph?} \label{sec:frucht}

The first theorem of this section answers the question from the title of the section.
The second theorem answers the same question for edge-colored digraphs; the 
conditions {in this second characterization} being a subset of the conditions for graphs.
{In case of graphs it is also observed that the partial automorphism monoid uniquely determines the graph
up to forming the complement.}

\begin{Thm}[Partial automorphism monoids of graphs]
\label{thm:pautgraph}
Given an inverse submonoid $\Sgrp$ of $\is{X}$, where $X$ is a finite set, 
$|X| \geq 2 $, there exists a graph 
$\Gamma=(X,E)$ whose partial automorphism monoid 
$\paut\Gamma$ is equal to 
$\Sgrp$ if and only if the following conditions hold:
\begin{enumerate} 
	\item \label{thm:pautgraphi}
	$\Sgrp$ is a full inverse submonoid of $\is{X}$,
	\item \label{thm:pautgraphii}
	for any compatible subset $A \subseteq \Sgrp$ of rank $1$ partial
	permutations, 
	if $\Sgrp$ contains the join of any two elements of $A$, then 
	$\Sgrp$ contains the join of the set $A$,
	\item \label{thm:pautgraphiii}
	the rank $2$ elements of $\Sgrp$ form at least one and at most two $\Drel$-classes,
	\item \label{thm:pautgraphiv}
	the $\Hrel$-classes of rank $2$ elements are nontrivial.
\end{enumerate}
\end{Thm}

\begin{Proof}
Suppose $\Gamma$ is a graph on the vertex set $X$ of size at least $2$, and take the inverse monoid $\Sgrp=\paut{\Gamma}$. We begin by proving that $\Sgrp$ has the asserted properties.
The simple observation that all partial identity permutations are partial automorphisms implies 
condition (\ref{thm:pautdigraphi}), while condition (\ref{thm:pautdigraphii}) follows from Proposition \ref{prop:joins1}.

For condition (\ref{thm:pautgraphiii}), recall that by Proposition \ref{prop:d}, the $\Drel$-classes correspond to the isomorphism classes of the induced subgraphs. 
If $\Gamma$ is {a graph on at least two vertices, then} any induced subgraph of $\Gamma$ with two vertices either has no edge or has a single edge between the two vertices. Since any two graphs on two vertices containing an edge are isomorphic, and so are any two graphs
on two vertices containing no edges, this yields at most two $\Drel$-classes of partial automorphisms
of rank $2$. Since $\Gamma$ contains at least one subgraph on at least two vertices,
it admits at least one $\Drel$-class of partial automorphisms of rank $2$.

{A} group $\Hrel$-clas{s o}f any $\Drel$-class {is} isomorphic to the automorphism group of 
{the induced subgraph on} the common domain. Both kinds of subgraphs with two vertices have a two-element automorphism group containing the identity map and a transposition; isomorphic to ${\Bbb Z}_2$. Hence, the $\Hrel$-classes of rank $2$ 
are non-trivial, proving condition (\ref{thm:pautgraphiv}).

Our proof of the converse statement is constructive: Given an inverse submonoid $\Sgrp$
of $\is{X}$ that satisfies the conditions listed in the theorem, we will construct a graph 
$\Gamma_\Sgrp$ on the set $X$ for which $\paut{\Gamma_\Sgrp} = 
\Sgrp$. 

We define $\Gamma_\Sgrp=(X,E)$ as follows. If $\Sgrp$ contains just one 
$\Drel$-class of rank $2$, denote it either $D_{e}$ or $D_{n}$, and 
if $\Sgrp$ contains two $\Drel$-classes of rank $2$, denote them by 
$D_{e}$ and $D_{n}$ (standing for `edges' and `non-edges', respectively; the choice of `which is which' can be made arbitrarily). 
For any two distinct vertices $v_1, v_2 \in X$, we le{t $\{v_1, v_2\} \in E$ i}f and only if the partial permutation 
$(v_1)\vee(v_2)$ belongs to $D_{e}$. 

Next, we show that $\Sgrp$ and $\paut{\Gamma_\Sgrp}$ coincide on their elements of rank at most $2$.
Let $\varphi \in \Sgrp$ be a permutation of rank at most $2${, and we intend to verify that 
$\varphi\in\paut{\Gamma_\Sgrp}$}. 
If $\varphi$ has rank $1$ or $0$, then it is a partial graph automorphism of 
$\Gamma_\Sgrp$ (and of any graph on $X$), as all {sub}graphs induced by one vertex are isomorphic, 
and the empty map is a partial graph automorphism of any graph by definition. 
Next, suppose $\varphi$ has rank $2$, that is, $\varphi=[v_1\, u_1)\vee[v_2\, u_2)$, 
where $u_1 \neq u_2$ (and therefore $v_1 \neq v_2$). 
Since, $\varphi\varphi^{-1}=(v_1) \vee (v_2)$ and $\varphi^{-1}\varphi=(u_1) \vee (u_2)$, 
$\varphi\varphi^{-1} \Drel \varphi^{-1}\varphi$ implies $(v_1) \vee (v_2) \Drel (u_1) \vee (u_2)$ 
in $\mathcal S$.
Thus, by the definition of 
$\Gamma_\Sgrp$, we obtain the equivalence $(u_1, u_2)\in E$ 
if and only if $(v_1, v_2)\in E$. 
Hence, $\varphi \in \paut{\Gamma_\Sgrp}$. 

Next we turn to proving that each partial automorphism of $\Gamma_\Sgrp$ of rank at most $2$ belongs to $\Sgrp$.
Since $\Sgrp$ is a full inverse submonoid of $\is{X}$, all idempotents of $\is{X}$ belong to $\Sgrp$.
In particular, the empty map belongs to $\Sgrp$, and it is, of course, the single rank $0$ partial permutation in $\Sgrp$.
We establish that every rank $1$ partial permutatio{n $[v_2v_1)$ of $X$ is} contained in $\Sgrp$. 
By condition (\ref{thm:pautgraphiv}), the $\Hrel$-class of the rank $2$ idempotent $(v_2)\vee(v_1)$ is nontrivial, therefore, 
besides the idempotent itself, it contains another partial permutation which is an automorphism of
the induced subgraph on $\{v_1,v_2\}$.
Necessarily, this partial permutation swaps the two vertices $v_1,v_2$.
This implies that {the cycle $(v_2 v_1)$ belongs} to $\Sgrp$
for every pair $v_1,v_2$ of distinct elements of $X$. 
Since any rank $1$ partial permutation arises as a restriction of such a partial permutation, the claim about rank $1$ partial permutations follows.

Now suppose $\varphi$ is a rank $2$ partial automorphism of $\Gamma_\Sgrp$, that is, 
$\varphi=[v_1\, u_1)\vee[v_2\, u_2)$ for some $u_1,u_2,v_1,v_2\in X$ with $u_1\not= u_2$ and $v_1\not= v_2$.
Then $u_1$ and $u_2$ are connected by an edge if and only if $v_1$ and $v_2$ are, whence $(v_1) \vee (v_2) \Drel (u_1) \vee (u_2)$ in $\Sgrp$ by the definition of $\Gamma_\Sgrp$. 
Therefore, there exists an element $\psi$ in $\Sgrp$ with domain $\{u_1, u_2\}$ and range $\{v_1, v_2\}$. 
Moreover, since the rank $2$ $\Hrel$-classes of $\Sgrp$ are nontrivial by condition (\ref{thm:pautgraphiv}), 
the $\Hrel$-class of $\psi$ contains both partial permutations wit{h d}omain $\{u_1,u_2\}$ and range $\{v_1,v_2\}$, 
in particular, $\varphi \in \Sgrp$.

We have shown that $\Sgrp$ and $\paut{\Gamma}$ coincide on their elements of rank at most $2$, and so, by applying Proposition \ref{prop:joins2}, we obtain $\paut{\Gamma_\Sgrp}=\Sgrp$.
\end{Proof}

As we have observed for any graph $\Gamma$, the idempotents of the two rank $2$ $\Drel$-classes of $\paut{\Gamma}$ correspond to pairs of vertices forming edges and non-edges. 
Consequently, the graph $\Gamma$ with $\Sgrp=\paut{\Gamma}$ is uniquely determined up to forming complement.
Hence the following can be deduced.

\begin{Cor}
\label{cor:paut-eq}
{If $\Gamma,\Gamma'$ are graphs with $\paut\Gamma=\paut{\Gamma'}$ then either $\Gamma=\Gamma'$ or
$\Gamma=\widetilde{\Gamma'}$.}
\end{Cor}

\begin{Thm}[Partial automorphism monoids of edge-colored digraphs]
\label{thm:pautdigraph}
Given an inverse submonoid $\Sgrp$ of $\is{X}$, where $X$ is a finite set, there exists an edge-colored digraph $\Gamma=(X,E_1,\ldots,E_l)$ whose partial automorphism 
monoid $\paut\Gamma$ is equal to $\Sgrp$ if and only if the following conditions hold:
\begin{enumerate}
	\item \label{thm:pautdigraphi}
	$\Sgrp$ is a full inverse submonoid of $\is{X}$,
	\item \label{thm:pautdigraphii}
	for any compatible subset $A \subseteq \Sgrp$ of rank $1$ partial permutations, 
	if $\Sgrp$ contains the join of any two elements of $A$, then 
	$\Sgrp$ contains the join of the set $A$.
\end{enumerate}
\end{Thm}

\begin{Proof}
Suppose $\Gamma=(X,E_1,\ldots,E_l)$ is an edge-colored digraph, and consider $\Sgrp=\paut{\Gamma}$. 
As before, the fact that all partial identical permutations are partial automorphisms implies 
that $\Sgrp$ satisfies condition (\ref{thm:pautdigraphi}), while condition (\ref{thm:pautdigraphii}) follows from Proposition \ref{prop:joins1}.

The converse part of our proof is again constructive.
Suppose $\Sgrp \subseteq \is X$ has the asserted properties, and let us define an edge-colored graph $\Gamma_\Sgrp=(X,E_c\ (c\in C))$ as follows. 
Let the color palette {be $C=C_1 \cup C_2$ with $C_1 \cap C_2=\emptyset$}, where $C_1$ indexes the set of rank $1$, $C_2$ the set of rank $2$ $\Drel$-classes of $\Sgrp$. 
For every $c\in C_1$, let $E_c$ consist of all loops $(v,v)$ where $(v)$ is in the $\Drel$-class $D_c$.
Moreover, for every $c\in C_2$, let us choose and fix vertices $v_1^{c}, v_2^{c}\in X$ with 
$(v_1^c)\vee(v_2^c) \in D_c$, and 
define $E_c$ to consist of all edges $(u_1,u_2)$ such that 
$[u_1\, v_1^c)\vee[u_2\, v_2^c) \in \Sgrp$.
Note that this implies $(v_1^c, v_2^c)\in E_c$ and 
$(v_1^c)\vee(v_2^c) \Drel (u_1)\vee(u_2)$ by definition. 
The converse of the latter relation is not quite true. 
However, if 
$(v_1^c)\vee(v_2^c) \Drel (u_1)\vee(u_2)$ in $\Sgrp$, then at least one of the partial permutations 
$[u_1\, v_1^c)\vee[u_2\, v_2^c)$ and 
$[u_2\, v_1^c)\vee[u_1\, v_2^c)$ is in $\Sgrp$, and
therefore, for any such $u_1, u_2$, at least one of the edges $(u_1, u_2)$ and $(u_2,u_1)$ belongs to $E_c$.
Both belong to $E_c$ if any only if the $\Hrel$-classes in $D_c$ are nontrivial.

To complete the proof, we again intend to show that $\Sgrp$ and $\paut{\Gamma_\Sgrp}$ coincide on their elements 
of rank at most $2$.
Suppose {first} that $\varphi \in \Sgrp$ is of rank $1$ or $2${, and check that 
$\varphi\in\paut{\Gamma_\Sgrp}$}. 
In the first case, $\varphi=[v_2\, v_1)$, and so $(v_1) \Drel (v_2)$ in $\Sgrp$. 
Therefore, by the definition of $\Gamma_\Sgrp$, {both} subgraphs induced {on $\{v_1\}$ and 
$\{v_2\}$ contain} a single loop of the same color, making $\varphi$ a partial automorphism. 
In the second case, suppose that 
$\varphi=[v_1\, u_1)\vee[v_2\, u_2)$ with $u_1 \neq u_2$ and $(u_1,u_2) \in E_c$. 
Then  
$[u_1\, v_1^c)\vee[u_2\, v_2^c) \in \Sgrp$, and 
$\big([v_1\, u_1)\vee[v_2\, u_2)\big)\big([u_1\, v_1^c)\vee[u_2\, v_2^c)\big)=
[v_1\, v_1^c)\vee[v_2\, v_2^c) \in \Sgrp$ follows from $u_1\neq u_2$.
This implies $(v_1, v_2)\in E_c$. 
Similar arguments interchanging the $u$'s and $v$'s imply that $\varphi$ is indeed in $\paut{\Gamma_\Sgrp}$. 

Now suppose that $\varphi \in \paut{\Gamma_\Sgrp}$ is of rank $1$ or $2$, we need to verify that $\varphi\in\Sgrp$. 
If $\varphi=[v_2\, v_1)$, then the subgraphs induced on $\{v_1\}$ and $\{v_2\}$ are isomorphic.
This implies by the definition of $\Gamma_\Sgrp$ that $(v_1) \Drel (v_2)$ in $\Sgrp$, and hence $\varphi \in \Sgrp$. 
Otherwise, let $\varphi=[v_1\, u_1)\vee[v_2\, u_2)$ with $u_1 \neq u_2$, and let $D_c$ be the $\Drel$-class of $\Sgrp$ containing $(v_1)\vee(v_2)$. 
As we have seen above, this implies $(u_1, u_2) \in E_c$ or $(u_2, u_1) \in E_c$.
If $(u_1, u_2) \in E_c$, that is, $[u_1\, v_1^c)\vee[u_2\, v_2^c) \in \Sgrp$, then $(v_1, v_2) \in E_c$ 
as well, 
therefore $[v_1\, v_1^c)\vee[v_2\, v_2^c) \in \Sgrp$, and 
$\big([v_1\, v_1^c)\vee[v_2\, v_2^c)\big)\big([u_1\, v_1^c)\vee[u_2\, v_2^c)\big)^{-1}=
[v_1\, u_1)\vee[v_2\, u_2) \in \Sgrp$ 
follows. 
I{f $(u_2, u_1) \in E_c$, then} $[v_1\, u_1)\vee[v_2\, u_2) \in \Sgrp$ can be similarly deduced.

Thus we have seen that $\Sgrp$ and $\paut{\Gamma}$ coincide on their elements of rank at most $2$, and by applying Proposition \ref{prop:joins2}, we conclude $\paut{\Gamma_\Sgrp}=\Sgrp$.
\end{Proof}

\begin{Rem}
Since any edge-colored digraph corresponds to a digraph admitting the same inverse monoid of partial automorphisms in which the colors are encoded with different numbers of edges,
partial automorphism monoids of (monochromatic) digraphs admit the same characterization as the one in Theorem \ref{thm:pautdigraph}. 
\end{Rem}

\section{When is an inverse monoid isomorphic to the partial automorphism monoid of a graph?}
\label{sec:frucht2}

{This section gives abstract characterizations for the inverse monoids of partial automorphism monoids of graphs, 
and more generally, of edge-colored digraphs.
It is also determined when the partial automorphism monoids of graphs are isomorphic.}

These characterizations are obtained from the theorems of the previous section. The link between the abstract inverse monoids and partial permutation monoids is established by a modified version of the so-called Munn-representation \cite[Chapter 5.2]{Lawson}.

As mentioned before, all inverse monoids can be represented by partial permutations. 
The standard representation is the Wagner--Preston representation, which represents 
$\Sgrp$ as an inverse submonoid of $\is{\Sgrp}$ faithfully, and is the analogue of the Cayley representation. 
Nevertheless, as in the case of groups, there also exist representations on smaller sets which are faithful for special classes of inverse monoids. In this section, we shall use
one of them which represents an inverse monoid on the set of its idempotents.  

Let $\Sgrp$ be an inverse monoid, and for any $e \in E(\Sgrp)$, let $[e]=\{f \in E(\Sgrp): f \leq e\}$, the \emph{order ideal} of $E(\Sgrp)$ generated by $e$. 
The \emph{Munn representation} of $\Sgrp$ is the map
$$\delta_\Sgrp \colon \Sgrp \to \is{E(\Sgrp)},\ s \mapsto \delta_s,$$ where
$$\delta_s \colon [s^{-1}s] \to [ss^{-1}], \ e \mapsto ses^{-1}.$$
An inverse monoid $\Sgrp$ is called \emph{fundamental} if its Munn representation is faithful (injective). 
It can be proven that each symmetric inverse monoid is fundamental, and each full inverse submonoid of a fundamental inverse semigroup is also fundamental. 
Hence the partial automorphism monoid $\paut \Gamma$ of any graph, digraph, or colored digraph $\Gamma$, 
is fundamental.

Let $\Sgrp$ be a finite inverse monoid, and let $X$ denote the set of atoms of $E(\Sgrp)$. 
The {\em restricted Munn representation} 
of $\Sgrp$ is the Munn representation of $\Sgrp$ restricted to the atoms of $E(\Sgrp)$, 
that is, it is the {homomorphism}
$$\alpha_\Sgrp \colon \Sgrp \to \is{X},\ s \mapsto \delta_s|_{X},$$ where
$$\delta_s|_{X} \colon [s^{-1}s] \cap X \to [ss^{-1}] \cap X,\ e \mapsto ses^{-1}.$$

Note that $\alpha_\Sgrp$ is well defined for any finite inverse monoid, as 
the atoms of an order ideal of $E(\Sgrp)$ are atoms of $E(\Sgrp)$, and
an order isomorphism between posets maps atoms to atoms. 
Naturally, this representation can only be faithful (injective) if $\Sgrp$ is fundamental. 
However, this is not a sufficient condition in general. 
The next lemma states that within the class of finite Boolean inverse monoids, the two properties are indeed equivalent.

\begin{Prop}
\label{prop:boolfund}
For a finite Boolean inverse monoid $\Sgrp$, the restricted Munn representation of $\Sgrp$ is 
{faithful (injective)}
if and only if $\Sgrp$ is fundamental.
\end{Prop}

\begin{Proof}
If $\Sgrp$ is a finite Boolean inverse monoid, then any principal order ideal of $E(\Sgrp)$ is a (finite) Boolean algebra, and an order isomorphism between order ideals is also a Boolean algebra isomorphism (see \cite{lattice}).
Since finite Boolean algebras are generated by their atoms, an isomorphism $\delta \colon I \to J$ between Boolean algebras is uniquely determined by its restriction to the atoms of $I$. 
Hence, for any elements $s, t \in \Sgrp$, we have $\delta_s|_{X}=\delta_t|_{X}$ if and only if $\delta_s=\delta_t$, making the restricted Munn representation faithful (injective) if and only if the Munn representation is.
\end{Proof}

There is another reason why the restricted Munn representation is a very natural representation to consider for the partial automorphism monoids, and this is the following statement.

\begin{Prop}
\label{prop:repr}
For a partial automorphism monoid $\Sgrp=\paut{\Gamma}$ of an {edge}-colored digraph $\Gamma$, the restricted Munn representation $\alpha_\Sgrp$ and the representation of $\Sgrp$ on $V(\Gamma)$ are essentially the same.
More precisely, 
if $X$ is the set of atoms of $E(\Sgrp)$,
then the map 
$$\xi\colon V(\Gamma)\to X,\ v\mapsto \id_{\{v\}}$$
is a bijection with the property that, for every $\varphi\in \Sgrp$, we have 
\begin{equation}
\label{prop:repr1}
{\xi(\dom\varphi)=\dom(\alpha_\Sgrp(\varphi)),\quad \xi(\ran\varphi)=\ran(\alpha_\Sgrp(\varphi)),}
\end{equation}
{and}
\begin{equation}
\label{prop:repr2}
\xi(\varphi(v))=(\alpha_\Sgrp(\varphi))(\xi(v))\quad {\hbox{for any}\  v\in V(\Gamma).}
\end{equation}
\end{Prop}

\begin{Proof}
The semilattice of idempotents $E(\Sgrp)$ of $\Sgrp$ consists of the identical maps on the subsets of $V(\Gamma)$, that is,
$\Theta\colon \mathcal P(V(\Gamma))\to E(\Sgrp),\ Y\mapsto \id_{Y}$
is an isomorphism from the {meet-semilattice} $(\mathcal{P}(V(\Gamma));\cap)$ of the Boolean algebra 
$\mathcal{P}(V(\Gamma))$
to $E(\Sgrp)$.
Since the atoms of $\mathcal{P}(V(\Gamma))$ are the singleton subsets, and the set of atoms of $E(\Sgrp)$ is
$X=\{\id_{\{v\}}: v\in V(\Gamma)\}$,
the restriction of $\Theta$ to the sets of atoms induces the bijection
$\xi\colon  V(\Gamma)\to X,\ v\mapsto \id_{\{v\}}$.
This provides a natural identification of the atoms of $E(\Sgrp)$ with the vertices of $\Gamma$. 

The equalities (\ref{prop:repr1}) and (\ref{prop:repr2}), which we prove in this paragraph, say that, under this identification, 
$\varphi$ and $\alpha_\Sgrp(\varphi)$ are the same for every $\varphi\in \Sgrp$.
By definition, $\dom(\alpha_\Sgrp(\varphi))=[\varphi^{-1}\varphi]\cap X=
\{\id_{\{v\}}: v\in \dom\varphi\}=\xi(\dom\varphi)$, and the second equality in (\ref{prop:repr1}) follows similarly.
Furthermore, if $v\in\dom\varphi$, then we have
$(\alpha_\Sgrp(\varphi))\left(\id_{\{v\}}\right)
=\delta_{\varphi}|_{X}\left(\id_{\{v\}}\right)
=\varphi\id_{\{v\}}\varphi^{-1}=\id_{\{\varphi(v)\}}$.
\end{Proof}

We are now ready to prove the main theorems of the section.

\begin{Thm}[Partial automorphism monoids of graphs, up to isomorphism]
\label{thm:pautgraph2}
Given a finite inverse monoid $\Sgrp$, there exists a finite graph whose partial automorphism monoid is isomorphic to $\Sgrp$ 
if and only if the following conditions hold:
\begin{enumerate}
	\item \label{thm:pautgraph2i}
	$\Sgrp$ is Boolean,
	\item \label{thm:pautgraph2ii}
	$\Sgrp$ is fundamental,
	\item \label{thm:pautgraph2iii}
	for any subset $A \subseteq \Sgrp$ of compatible $0$-minimal elements, if all $2$-element subsets of $A$ have {joins in} $\Sgrp$, then the set $A$ has a join in $\Sgrp$,
	\item \label{thm:pautgraph2iv}
	$\Sgrp$ has at most two $\mathcal D$-classes of height $2$,
	\item \label{thm:pautgraph2v}
	the $\mathcal H$-classes of the height $2$ $\mathcal D$-classes of $\Sgrp$ are nontrivial.
\end{enumerate}
\end{Thm}

\begin{Proof}
For any graph $\Gamma$, we have seen 
that $\paut{\Gamma}$ is Boolean and fundamental, proving condition{s (\ref{thm:pautgraph2i}) 
and (\ref{thm:pautgraph2ii})}.

Since the $0$-minimum elements of $\paut{\Gamma}$ are just the rank $1$ elements, and the $\Drel$-classes of height $2$ are just the rank $2$ $\Drel$-classes, conditions (\ref{thm:pautgraph2iii})--(\ref{thm:pautgraph2v}) 
are immediate consequences of our previous observations on partial automorphism monoids; see properties (\ref{thm:pautgraphii})--(\ref{thm:pautgraphiv}) in Theorem \ref{thm:pautgraph}.

Conversely, let $\Sgrp$ be an inverse monoid having the asserted properties, and let $X$ denote the set of all atoms of $E(\Sgrp)$. 
By propert{ies} (\ref{thm:pautgraph2i}) and (\ref{thm:pautgraph2ii}), Proposition \ref{prop:boolfund} implies 
that the restricted Munn representation $\alpha_\Sgrp$ of $\Sgrp$ embeds $\Sgrp$ into $\is{X}$. 
We check that the inverse submonoid $\alpha_\Sgrp(\Sgrp)$ of $\is{X}$ satisfies conditions
{(\ref{thm:pautgraphi})--(\ref{thm:pautgraphii})} of Theorem \ref{thm:pautgraph}. 
Since the rank $1$ elements of $\alpha_\Sgrp(\Sgrp)$ are clearly $0$-minimal, 
condition (\ref{thm:pautgraphii}) is immediate from property (\ref{thm:pautgraph2iii}) on joins in $\Sgrp$. 
To verify condition (\ref{thm:pautgraphi}), that is, that $\alpha_\Sgrp(\Sgrp)$ is a full inverse 
submonoid of $\is{X}$, note that for any element $e \in X$, 
$\id_{\{e\}}=\alpha_\Sgrp(e) \in \alpha_\Sgrp(\Sgrp)$.
Since $E(\Sgrp)$ is a meet-semilattic{e o}f a Boolean algebra with $|X|$ atoms, the same holds also for $E(\alpha_\Sgrp(\Sgrp))$.
However, $E(\is{X})$ is also a Boolean algebra of the same size which contains $E(\alpha_\Sgrp(\Sgrp))$.
Hence $E(\alpha_\Sgrp(\Sgrp))=E(\is{X})$, and $\alpha_\Sgrp(\Sgrp)$ is, indeed, a full inverse submonoid of $\is{X}$. 
Since $\alpha_\Sgrp(\Sgrp)$ is a full inverse submonoid of $\is{X}$, its $\Drel$-classes of height $2$ are just
the rank $2$ $\Drel$-classes, therefore properties (\ref{thm:pautgraph2iv})--(\ref{thm:pautgraph2v})
imply that $\alpha_\Sgrp(\Sgrp)$
satisfies conditions (\ref{thm:pautgraphiii})--(\ref{thm:pautgraphiv}) of Theorem \ref{thm:pautgraph} as well.

Theorem \ref{thm:pautgraph} therefore shows the existence of a required $\Gamma$, thus completing
the proof.
\end{Proof}

Using exactly the same arguments and Theorem \ref{thm:pautdigraph}, one obtains the following theorem for edge-colored digraphs. 

\begin{Thm}[Partial automorphism monoids of edge-colored digraphs, up to isomorphism]
\label{thm:pautdigraph2}
Given a finite inverse monoid $\Sgrp$, there exists a finite edge-colored digraph whose partial automorphism monoid 
is isomorphic to $\Sgrp$ if and only if the following conditions hold:
\begin{enumerate}
	\item 
	$\Sgrp$ is Boolean,
	\item 
	$\Sgrp$ is fundamental,
	\item 
	for any subset $A \subseteq \Sgrp$ of compatible $0$-minimal elements, if all $2$-element subsets of $A$ have {joins in} $\Sgrp$, then the set $A$ has a join in $\Sgrp$.
\end{enumerate}
\end{Thm}

{We can combine our results to determine when the partial automorphism monoids of graphs are isomorphic.}

\begin{Cor}
{For any graphs $\Gamma$ and $\Gamma'$, the partial automorphism monoids $\paut\Gamma$ and $\paut{\Gamma'}$ are isomorphic if and only if either
$\Gamma$ and $\Gamma'$, or $\Gamma$ and $\widetilde{\Gamma'}$ are isomorphic.}
\end{Cor}

\begin{Proof}
The `if' part of the statement is obvious.
To see the converse, suppose that $\iota\colon \paut{\Gamma}\to \paut{\Gamma'}$ is an isomorphism, and put
$X=V(\Gamma),\ X'=V(\Gamma')$.
{The sets of elements of rank at most $1$ form inverse subsemigroups $S_1$ and $S'_1$ in $\paut{\Gamma}$ and
$\paut{\Gamma'}$, respectively.
The isomorphism $\iota$ maps idempotents to idempotents and preserves their heights,
and so it also preserves ranks of elements.
Therefore $\iota(S_1)=S'_1$, and
$\iota$ restricts to a bijection between the sets of 
rank $1$ idempotents of $\paut{\Gamma}$ and $\paut{\Gamma'}$.
Combining this bijection with the bijections $v\mapsto \id_{\{v\}}$ from the sets $X$ and $X'$ to the latter sets,
we obtain a bijection $\overline{\iota}\colon X\to X'$,
where $\id_{\{\overline{\iota}(v)\}}=\iota(\id_{\{v\}})$ for every $v\in X$.
Recall that $S_1\setminus \{\id_{\emptyset}\}$ and $S'_1\setminus \{\id_{\emptyset}\}$ constitute $\Drel$-classes in $\is{X}$ and $\is{X'}$, respectively.
Since isomorphisms map $\Lrel$-related ($\Rrel$-related) elements to elements with the same property, and the
$\Hrel$-classes of $S_1$ and $S'_1$ are singletons, we obtain that
$\iota([v\,u))=[\overline{\iota}(v)\,\overline{\iota}(u))$
for any $u,v\in X$.}

Now let us define a graph $\Gamma''=(X,E'')$ in the way that $\{u,v\}\in E''$ for distinct $u,v\in X$ if and only if
$\{\overline{\iota}(u),\overline{\iota}(v)\}\in E'$.
It is clear that $\Gamma''$ is isomorphic to $\Gamma'$.
We intend to verify that the rank $2$ elements of $\paut{\Gamma''}$ and  $\paut{\Gamma}$ coincide.

Let $\varphi=[v_1\,u_1)\vee [v_2\,u_2)$ be a rank $2$ element of $\is{X}$.
By definition, $\varphi\in\paut{\Gamma''}$ if and only if either 
both
$\{\overline{\iota}(u_1),\overline{\iota}(u_2)\}$ and $\{\overline{\iota}(v_1),\overline{\iota}(v_2)\}$ are edges, or
both are non-edges in $\Gamma'$, 
that is, if and only if 
$[\overline{\iota}(v_1)\,\overline{\iota}(u_1))\vee [\overline{\iota}(v_2)\,\overline{\iota}(u_2))\in\paut{\Gamma'}$.
As isomorphisms are easily seen to respect joins, we obtain that 
$$[\overline{\iota}(v_1)\,\overline{\iota}(u_1))\vee [\overline{\iota}(v_2)\,\overline{\iota}(u_2))=
\iota\big([v_1\,u_1)\big)\vee \iota\big([v_2\,u_2)\big)=
\iota\big([v_1\,u_1)\vee [v_2\,u_2)\big)=\iota(\varphi).$$
Hence $\varphi\in\paut{\Gamma''}$ if and only if $\varphi\in \iota^{-1}\big(\paut{\Gamma'}\big)=\paut{\Gamma}$,
and so the rank $2$ elements of $\paut{\Gamma''}$ and $\paut{\Gamma}$, indeed, coincide.

Since both in $\paut{\Gamma''}$ and $\paut{\Gamma}$, the elements of rank at most $1$ are just the elements of the same kind in $\is X$, we deduce by Proposition \ref{prop:joins2} that $\paut{\Gamma}=\paut{\Gamma''}$.
Thus Corollary \ref{cor:paut-eq} implies that either $\Gamma=\Gamma''$ or $\Gamma=\widetilde{\Gamma''}$, and
since $\Gamma''$ is isomorphic to $\Gamma'$, our corollary follows.
\end{Proof}

Motivated by this result and the Graph Reconstruction Conjecture, we propose a related question.

\begin{Probl}
Given graphs $\Gamma_1$ and $\Gamma_2$, if their decks $\operatorname{Deck}(\Gamma_1)$ and $\operatorname{Deck}(\Gamma_2)$ coincide up to isomorphism and complementation, 
do $\Gamma_1$ and $\Gamma_2$ necessarily coincide up to isomorphism and complementation?
\end{Probl}

This can be translated to the language of semigroups in the following way. Note that given a graph $\Gamma$, the partial automorphism monoids $\paut{\Gamma-v_i}$ ($v_i \in V$) are inverse submonoids of $\paut{\Gamma}$, specifically, the maximal Boolean inverse submonoids of height $|V(\Gamma)|-1$. Let $\operatorname{Deck}(\paut{\Gamma})$ denote the multiset of such inverse submonoids of $\paut{\Gamma}$.

\begin{Probl}
\label{probl:last}
Given two inverse monoids $\paut{\Gamma_1}$ and $\paut{\Gamma_2}$, if $\operatorname{Deck}(\paut{\Gamma_1})$ and $\operatorname{Deck}(\paut{\Gamma_2})$ coincide up to isomorphism, do we have $\paut{\Gamma_1} \cong \paut{\Gamma_2}$?
\end{Probl}

We note that the above problems can fail for small graphs. E.g., take the following two graphs on four vertices:

\begin{center}
\begin{figure}[H]
\includegraphics[width=0.35\linewidth]{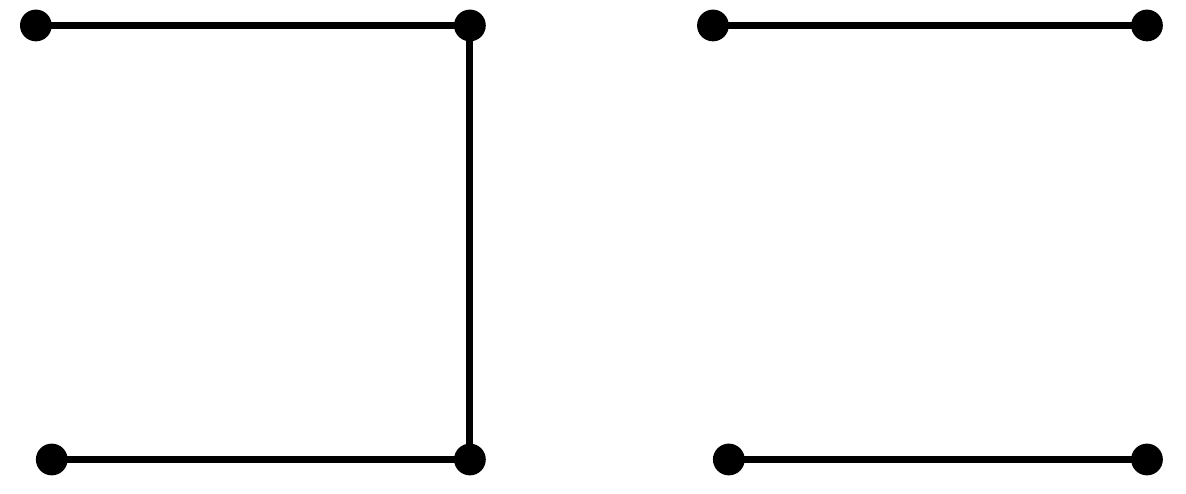}
\end{figure}
\end{center}

One can check that $\operatorname{Deck}(\paut{\Gamma_1})$ and $\operatorname{Deck}(\paut{\Gamma_2})$ are the same, but $\paut{\Gamma_1}$ and $\paut{\Gamma_2}$ are different.
However, we believe that for large enough graphs, the answer to Problem \ref{probl:last} may be positive.

\end{document}